\documentclass {article}
 \advance \textheight by \topskip
 \newcommand { \be} {\begin {equation}} 
 \newcommand { \en} { \end {equation}} 
 \newcommand { \bea} {\begin {eqnarray}} 
 \newcommand { \ena} { \end {eqnarray}} 
 \newcommand { \bean} {\begin {eqnarray*}} 
 \newcommand { \enan} { \end {eqnarray*}} 
 \newcommand { \bfig} {\begin {figure}} 
 \newcommand { \efig} { \end {figure}} 
 \newcommand { \Proof} {\smallskip \noindent {\bf Proof: }}
 \newtheorem {theorem} {Theorem} 
 \newtheorem {proposition} {Proposition} 
 \newtheorem {lemma} [proposition] {Lemma}
 
 \newtheorem {corollary} [proposition] {Corollary} 
 \newcommand { \Remark} {\smallskip \noindent {\bf Remark: }} 
 \newcommand { \qed} { $ \hfill \clubsuit $ \smallskip}
 \def \ack { \noindent {\bf Acknowledgments \par   }} 
\newdimen\figcenter
\newdimen\texpscorrection
\def\putfigurewithtex#1 #2 #3 #4
{
{
\figcenter=\hsize\relax
\advance\figcenter by -#4truecm
\divide\figcenter by 2
\vskip #3truecm
\noindent\hskip\figcenter
\includegraphics{#1}
{\hskip\texpscorrection\input #2 }
}
}     
\begin {document}
\title { Elliptic Islands on the Elliptical Stadium} 
\author 
{Sylvie Oliffson Kamphorst 
\hskip 0.6 truecm 
S\^onia Pinto de Carvalho 
\\ 
Departamento de Matem\'atica, ICEx, UFMG \\ 
C.P. 702, 30161--970, Belo Horizonte, Brasil.} 
\date {} 
\maketitle 
\footnotetext{AMS subject classification: 37E40, 70K42}
\begin{abstract}
We investigate the existence of elliptic islands for a special family of
periodic orbits of a two-parameter family of maps $T_{a,h}$, corresponding to
the billiard problem on the elliptical stadium.

The hyperbolic character of those orbits were studied in \cite{kn:cana} for
$1<a<\sqrt 2$ and here we look for the elliptical character for every $a>1$.

We prove that, for $a<\sqrt 2$, the lower bound for chaos $h=H(a)$ found in
\cite{kn:cana} is the upper bound of ellipticity for this special family.

For $a>\sqrt 2$ we prove that there is no upper bound on $h$ for the existence
of elliptic islands.

The main results we use are Birkhoff Normal Form and Moser's Twist Theorem.
 \end{abstract}

\section
{Introduction} 

The elliptical stadium is a plane region bounded by a curve $ \Gamma $, 
constructed by joining two half-ellipses, with major axes $ a>1 $ and minor axes
$ b = 1 $, by two straight segments of equal length $ 2h $
(see fig. \ref {fig:estadio}). 
\bfig [h]
\putfigurewithtex bilhar.ps bilhar.txt 4 16 
\caption {The elliptical stadium.} 
\label {fig:estadio} 
\efig 

The billiard on the elliptical stadium consists in the study of the free 
motion of a point particle inside the stadium, being reflected elastically 
at the impacts with $ \Gamma $. 
Since the motion is free inside $ \Gamma $, it is determined 
either by two consecutive points of reflection at $ \Gamma $ or by the point of 
reflection and the direction of motion immediately after each collision. 

For fixed $ a $ and $ h $, let $ s \in [0,L) $ be the arc length parameter for
$ \Gamma $ 
and the direction of motion be given by the angle $ \beta $ with the normal to
the boundary at the impact point. The billiard defines an invertible map
$ T_{a, h} $ from the annulus $ { \cal A} = [0,L) \times (- \pi / 2, \pi / 2) $
into itself, preserving the measure $ d \mu = \cos \beta \, d \beta \, ds $. 
Since $ \Gamma $ is globally $ C^1 $ but not $ C^2 $,
$ T_{a, h} $ is a homeomorphism (see, for instance, \cite {kn:mar}) and if 
$ (s_0, \beta_0) $ and $ (s_1, \beta_1) =   T_{a, h} (s_0, \beta_0)
\in { \cal A} $ such that $ \Gamma $ is analytic in some
neighborhoods of $ s_0 $ and $ s_1 $, then clearly $ T_{a, h} $ is analytic in
some neighborhoods of $ (s_0, \beta_0) $ and $ (s_1, \beta_1) $.

For each $ (a, h) $, $ ({ \cal A}, \mu, T_{a, h}) $
defines a discrete dynamical system, whose dynamics depends on the values of
$ a $ and $ h $. For instance, when $ h = 0 $ we have an ellipse and the
billiard is integrable.

When $ h \neq 0 $, two main features appear. If $ a< \sqrt 2 $, Donnay \cite
{kn:don}  proved  that the billiard on the elliptical stadium is chaotic (in
the sense of non-vanishing Lyapunov exponents) when $h$ is sufficiently large.
Lower bounds for $h $ for this behaviour were found by Markarian and
ourselves in \cite{kn:cmp} and by Canale, Markarian and ourselves in
\cite{kn:cana}. In the present work we show that the lower bound found in 
\cite{kn:cana} is optimal, in the sense that bellow it we can assure the
existence of elliptic islands of positive measure. 

In \cite{kn:bun} Bunimovich conjectured the existence of a stable periodic
orbit, with island of positive measure, for billiards such 
as the elliptical stadium with $ a> \sqrt 2 $ and $ h \neq 0 $.
In this work we
make some progress in this direction, proving that there is no upper bound
on $ h $ for the existence of elliptic islands if $ a> \sqrt 2 $. So, there
is no way to destroy the elliptic islands by just increasing the distance
between the half-ellipses.

To prove the existence of elliptic islands we extend the results about a
special family of periodic orbits, called pantographic, studied in \cite
{kn:cana} and find regions on the parameter
plane where at least one of its members is elliptic and stable,
so, with islands of positive measure in phase space.

\section {Pantographic orbits: existence and ellipticity} 
 
In this section we define the special family of periodic orbits and investigate
the existence and ellipticity of its members. This family has already been
investigated in \cite {kn:cana} for $a < \sqrt {2}$. Here we extend
that work for all $a >1 $. We will skip most of the proofs which can be found
in the work cited above. 

Given $a$, $h$ and a positive integer $n$, an $(n,a,h)$-pantographic orbit, 
denoted by $Pan(n,a,h)$, is a symmetric $(4+2n)$-periodic orbit, with exactly
2 impacts at each half-ellipse, joined by a vertical path, and crossing any
vertical line only twice. One example can be seen in figure~\ref{fig:pan}. 

\bfig [h] 
\putfigurewithtex pan.ps  pan.txt 4 18
\caption {The 10-periodic pantographic orbit ($n=3$).}  
\label {fig:pan} 
\efig 

Let the right half-ellipse of the stadium be parametrized by
$(x,y) = (a\cos \lambda +h,\sin \lambda)$ and $P$ be the point marked on 
figure~\ref{fig:pan}. Using the obvious symmetries (see figure~\ref
{fig:equacao}), the parameter $\lambda$ of $P$ must satisfy:
\be
\tan 2 \beta = 
\frac { a \tan \lambda} { a^2 \tan^2 \lambda - 1} =
\frac {h + a \cos \lambda} {n + \sin \lambda}
\,\,\, {\mbox {and}}\,\,\,  \tan \beta = \frac{\cos \lambda}{a\sin \lambda}
\label {eq:pan}
\en
where $ \beta > 0 $ is, as defined above, the angle of the trajectory from
$ P $, with the normal to the boundary.

\bfig [h] 
\putfigurewithtex coord.ps  coord.txt 8.5 16
\caption {} 
\label {fig:equacao} 
\efig
 
The following proposition gives the
region of existence of those pantographic orbits in the parameter plane 

\begin {proposition} 
$ $ \\  \vskip -1truecm
\begin {itemize}
\vskip -2truecm
\item 
$ Pan (0, a, h) $ and $ Pan (1, a, h) $ exist for every $ a > 1 $ and $ h > 0 $. 
\item 
For $ n \geq 2 $, $ Pan (n, a, h) $ exists for every $ 1 < a \leq 2 $ and
$ h > 0 $ or for every $ a>2 $ and $ h> (n-1) \sqrt {a (a-2)} $.
\label{prop:existence}
\end {itemize}
\end {proposition} 
  
\Proof
Equation~(\ref {eq:pan}) can be written as
\be 
n = \frac {a^2 t^2-1} {2a t} \, h + \frac { (a^2-2) t^2-1} {2t \sqrt {1 + t^2}}.
\label {eq:y} 
\en
where $ t = \tan \lambda $.
As proven in \cite {kn:cana}, for each integer $ n \geq 0 $,
this equation has a unique solution
$ t (n, a, h) = \tan \lambda (n, a, h) > 0 $ for every 
$ a>1 $ and $ h>0 $ and $ t (n, a, h) \in (\frac {1} {a}, + \infty) $.

For $ n = 0 $ and $ n = 1 $ this results implies the existence of the 
corresponding $Pan(n,a,h)$ for every $ a>1 $ and $ h>0 $.

However, for $ n \geq 2 $ one must also ask that 
$$
\tan 2 \beta \geq \frac {a \cos \lambda} {1 + \sin \lambda}
$$ 
in order to guarantee that the next impact point from $P$ is on the straight
part of the boundary. 
This is  equivalent to  
\be
0 \leq \sin \lambda \leq \frac {1} {a-1}
\label{eq:seno} 
\en
which is always true if $ 1<a \leq 2 $. 

If $ a>2 $ we rewrite (\ref{eq:seno}) 
as $ 0 \leq t = \tan \lambda \leq \frac {1} { \sqrt {a (a-2)}} $.
Since $ \frac { \partial t} { \partial h} (n, a, h)<0 $
(which is easily verified from (\ref{eq:y})), and 
$ \frac {1} {a}< \frac {1} { \sqrt {a (a-2)}} $, there
exists a unique $ \overline h $ such that 
 $ t (n, a, \overline h) = \frac {1} { \sqrt {a (a-2)}} $ and if
 $ h> \overline h $,
 $ \frac {1} {a} < t (n, a, h)< \frac {1} { \sqrt {a (a-2)}} $.
From (\ref {eq:y}), $ \overline h = (n-1) \sqrt {a (a-2)} $. 
\qed
 
For each fixed $ n $, we denote by  $ U_n $ the open region in the parameter
plane where $ Pan (n, a, h) $ exists, according to Proposition~\ref{prop:existence}.
For $ (a, h) \in U_n $, let $ s $  be the 
arc length corresponding to the point $P$ of $ Pan (n, a, h) $ and $ \beta $
the angle with the normal of the trajectory at this point as before.
Then $ T_{a, h}^{4 + 2n} \, (s, \beta) = (s, \beta) $ and the 
ellipticity of this orbit is determined by the eigenvalues of 
 $ { \cal D} T_{a, h}^{4 + 2n} | _{ (s, \beta)} $. 

As shown in \cite {kn:cana}, we can write 
$ { \cal D} T_{a, h}^{4 + 2n} | _{ (s, \beta)} = (M_1 M_2)^2 $ with 
$$ 
M_j = { 1 \over \cos \beta } 
\left (
\begin {array} {cc} 
l_j \, K - \cos \beta & l_j \\ 
K \, (l_j \, K - 2 \cos \beta) & l_j \, K -  \cos \beta 
\end {array} \right)
$$ 
and where 
$ l_1 $ is the length of the trajectory between two impacts with the same
half-ellipse, $ l_2 $ is the length of the trajectory between two impacts with
the different half-ellipses and $ K $ is the curvature of the ellipse at $ s $. 

Let 
$$ 
\Delta_n (a, h) = \left ({l_1 \, K \over \cos \beta} -1 \right ) \, 
 \left ({l_2 \, K \over \cos \beta} -1 \right ). 
$$ 
Since det $ (M_1 M_2) = 1 $ and tr$(M_1M_2)=4\Delta_n(a,h)-2$, it follows
that if $ 0< \Delta_n (a, h)< \frac {1} {2}< \Delta_n (a, h)<1 $  then $ Pan
(n, a, h) $ is { \em elliptic} (meaning that the eigenvalues of   $ { \cal D}
T^{4 + 2n} | _{ (s, \beta)} $ are unitary with non zero imaginary part).

The following lemma summarizes some properties of $ \Delta_n (a, h) $
and its technical proof has been postponed to the appendix.  

\begin {lemma}
For every $ n \geq 0 $, let 
$$ \tilde U_n = \{ (a, h) | \,1 < a < \sqrt 2, h>0 \}
\cup
\{ (a, h) | \, a \geq \sqrt 2, h > n a \sqrt {a^2-2} \equiv h^0_n(a) \}
\subset U_n \ . $$
The function $ \Delta_n (a, h) $ has the following properties: 
\begin {enumerate} 
\item
$ \Delta_n (a, h) | _{ \tilde U_n } > 0 $ 
\item
$ \frac { \partial \Delta_n} { \partial h} | _{ \tilde U_n } > 0 $ 
\item
$ {\displaystyle \lim_{h \rightarrow + \infty}} \Delta_n (a, h) = + \infty $ 
\item 
for $ 1 < a < \sqrt 2 $,
$ {\displaystyle \lim_{h \rightarrow 0}} \Delta_n (a, h) = 
L_n (a) = (\frac {2} {a^2}-1) (\frac {2 (n + 1)} {a^2}-1)> 0 $ \\
for $ a \geq \sqrt 2 $,
$ {\displaystyle \lim_{h \rightarrow n a \sqrt {a^2-2}} }\Delta_n (a, h) = 0 $.
\end {enumerate} 
\end {lemma}

For each $ 0 < c \leq 1 + 2n $, let $ \alpha_n^c $ be the unique solution of
$ L_n (a) = c $. We have that $ 1< \alpha_n^c < \sqrt 2 $ and if $ c < d $ then
$ \alpha_n^d < \alpha_n^c $.

It follows from the lemma that every level curve $ \Delta_n = c $ is given by a
graph $ h_n^c: (\alpha_n^c, + \infty) \to I \hspace {-.12truecm} R $ such that
$ \Delta_n (a, h_n^c (a)) = c $ and $ \Delta_n (a, h) < c $ if
$ h < h_n^c (a) $ and $ \Delta_n (a, h) > c $ if $ h > h_n^c (a) $.

The characterization of the region of ellipticity is then given by:
\begin {proposition}
For each fixed $ n $, the region in the parameter plane where $ Pan (n, a, h)
$ is elliptic is the union of two open adjacent strips, bounded by the graphs of 
$ h_n^0 (a), h_n^{1/2} (a), h_n^1 (a) $ and by the segments
$ \{ (a,0) | \, \alpha_n^1 < a < \alpha_n^{1/2} \} $ and 
$ \{ (a,0) |\, \alpha_n^{1/2} < a < \sqrt 2 \} $.
(see figure~\ref {fig:juntas}).
\label {prop:region}
\end {proposition}

\bfig
\putfigurewithtex elres.ps  elres.txt 7.5 16
\caption {
Region of ellipticity with lines of resonance up to order 4 for $ n =3 $ } 
\label {fig:juntas} 
\efig

For each $ n $ and $ (a, h) $ in the region of ellipticity of $ Pan (n, a, h) $,
let $ \mu $ and $ \overline \mu $ be the eigenvalues of
$ DT_{a, h}^{4 + 2n} | _{ (s, \beta)} $.
$ Pan (n, a, h) $ has a resonance of order $ k $ if $ \mu^k = 1 $, i.e., 
$ \mu = e^{i \frac {2j \pi} {k}}, j = 1,...,k-1 $.
It is easy to show that this corresponds, for $k>1$, to
$$
\Delta_n (a, h) 
= \frac{1}{2} \left ( {1 + \cos \frac {j \pi} {k} } \right )    \equiv c_{j k}.
$$ 
Clearly $ 0<c_{j k}<1 $ and we have the curves of resonance given by the graphs
of
$ h_n^{c_{j k}}: (\alpha_n^{c_{j k}}, + \infty) \to I \hspace {-.12truecm} R $.

We can then characterize the region of ellipticity with no resonances up to
order $ k $ by:
\begin {proposition}
\label{prop:reso}
For a fixed $ n $, the region in the parameter plane
where $ Pan (n, a, h) $ is elliptic with no resonances up to order $ k $ is a
finite union of open disjoint adjacent strips contained in the region of
ellipticity. (see figure~\ref {fig:juntas})
\end {proposition}

\section {Existence of elliptic islands}

 In order to establish if the elliptic periodic orbits described in the previous
section have invariant curves surrounding them, we will invoke the two 
classical results:

{ \bf Birkhoff Normal Form}: 
{ \em 
Let $ f $ be an area preserving map in $ C^l $ ($ l \geq 4 $ ) with a fixed
point  at the origin, with eigenvalues $ \mu $ and $ \overline { \mu} $, $ |
\mu | = 1 $. If for some integer $ q $ in $ 4 \leq q \leq l + 1 $ one has $
\mu^k \neq 1 $  for $ k = 1,2,...,q $ then there exists a real analytic
transformation taking $ f $ into the normal form
\be
\zeta \to f (\zeta, \overline \zeta) = 
\mu \zeta e^{i \tau (\zeta \overline \zeta)} + 
g (\zeta, \overline \zeta)
\label{eq:bnf}
\en
where $ \tau (\zeta \overline \zeta) =
\tau_1 | \zeta | ^2 + ... + \tau_s | \zeta | ^{2s} $, with
$ s = \left[ \frac {q} {2} \right] -1 $, is a real polynomial in
$ | \zeta | ^2 $ and $ g $ vanishes with its derivatives up to order $ q-1 $  at
$ \zeta = \overline \zeta = 0 $.}

{ \bf Theorem} ( {\bf Moser}, \cite{kn:mos}, p.56)
{ \em
If the polynomial $ \tau (| \zeta | ^2) $ does not vanishes identically,
$ \zeta = 0 $ is a stable fixed point}\,\, (which means that there are
invariant curves surrounding it, and so an elliptic island of positive
measure). 

For each fixed period $ n_0 $, we will then investigate the resonances of 
$ T_{a, h}^{4 + 2n_0} $ and the zeros of the coefficients
of its Birkoff normal form, near the Pantographic orbit.

Let us fix a period $n_0$ and a major axis
$ a_0 > \alpha_{n_0}^1 = \sqrt { (2 + 2n_0)/ (2 + n_0) } $. According to
proposition \ref{prop:reso}, $ Pan (n_0, a_0, h) $ is elliptic with
no resonances up to order $q$ if $ h $ is in the finite union of open
disjoint adjacent intervals denoted $ \cup I_j^q $.

Let $ \lambda (h), s (h), \beta (h) $ be respectively the parameter,
the arc length and
the angle with the normal for the point $ P $ of $ Pan (n_0, a_0, h) $.

\begin {lemma}
For any $q\geq 4$, $ \lambda (h), s (h) $ and $ \beta (h) $ are analytic
functions of $ h $ on each $ I_j^q $.
\label {lem:ana}
\end {lemma}
\Proof
The point $ P $ will belong to $ Pan (n_0, a_0, h) $ if $ t = \tan \lambda $
satisfies equation (\ref{eq:y}):
$$
n_0 = \frac {a_0^2 t^2-1} {2a_0 t} \, h +
\frac { (a_0^2-2) t^2-1} {2t \sqrt {1 + t^2}}.
$$ 
$ A (t) = \frac {a_0^2 t^2-1} {2a_0 t} $ and
$ B (t) = \frac { (a_0^2-2) t^2-1} {2t \sqrt {1 + t^2}} $ are analytical
functions of $ t $ and $ A (t) \neq 0 $, since $ \frac {1} {a_0} < t < \infty $.
So $ h = h (t) = \frac {n_0-B (t)} {A (t)} $ is analytic.
As this equation has a unique solution
for each $ h $, the inverse $ t = t (h) $  exists and is then
locally analytic for every $ h \in I_j^q $.

The functions $ \lambda (h) = \arctan \,t (h) $ and the corresponding arc length
of the ellipse $ s (h) = s (\lambda (h)) $ are then analytic. 
The same is true for $ \beta (h) = \beta (\lambda (h)) $.   \qed

For each fixed $h$, let $f$ be the translation of $T_{a_0, h}^{4 + 2n_0}$ by 
$ (s (h), \beta (h)) $.
The map $ f $ is clearly area preserving and analytic
in $ (s, \beta) $ on a neighbourhood of the origin. 
The eigenvalues of $ Df_{ (0,0)} $ are 
the same as those of $ DT_{a_0, h}^{4 + 2n_0} { (s (h), \beta (h))} $.
If $ h \in I_j^q $, $ f $ can be written in the Birkhoff normal form 
(\ref{eq:bnf}).
If one of the Birkkoff coefficients is not zero, $ f $ has an elliptic island surrounding $ (0,0)
$ and, by translation, there is an elliptic island surrounding $ Pan (n_0,
a_0, h) $.

\begin {lemma}
For $ 1 \le m \le \left [ \frac{q}{2} \right ] -1 $, $ \tau_m (h) $,
the $m$-th Birkhoff coefficient of $f$, is an analytic function of $ h $ on each
$ I_j^q $.
\end {lemma}
 \Proof
For each fixed $ h $, $ \tau_m (h) $ is a combination of the coefficients of the
($q-1$)-th jet of $f$ at $(0,0)$ (see, for instance,
 \cite {kn:moe}). The steps giving rise to the calculation of
 $ \tau_m (h) $ (complexification of the space, elimination of unwanted
terms) are analytical.  
So, if the coefficients of the jet $ J_{q-1} f_{ (0,0)} $ are
analytical in $ h $, then $ \tau_m (h) $ will also be an analytic function of $ h $.

Now these coefficients are combinations of the entries of $ DT_{a_0, h}^{4 + 2n_0} $
and their derivatives up to $(q-1)$-th order with respect to $ s $ and $ \beta
$, calculated at $ (s (h), \beta (h)) $.

Let $ (s, \beta) $ and $ (s', \beta') $ be two consecutive
impacts of a trajectory with the two different half-ellipses (with $ l \ge 0 $
impacts with the straight parts between them), or two consecutive impacts of a
trajectory with the same half-ellipse (with $ l = 0 $ ).
Then $ DT_{a_0, h}^{4 + 2n_0} (s, \beta) $, for 
$ (s, \beta) $ near to $ (s (h), \beta (h)) $ is a finite product of matrices of
the form (see, for instance, \cite {kn:mar})
$$ 
 { (-1)^l \over \cos \beta'} 
 \left (
\begin {array} {cc} 
L K - \cos \beta & L \\ 
K K' L - K' \cos \beta - K \cos \beta' & 
L K' - 
 \cos \beta' 
 \end {array} \right) 
$$ 
where $ K $ stands for the curvature of the ellipse at the impact and
$ L $  is the total length of the trajectory between the two impacts
with the half-ellipses. 

Since $ \cos \beta \neq 0 $ for $ \beta $ near $ \beta (h) $, all the entries
of the matrix above , as well as their derivatives up to
any order in $ s $ and $ \beta $, are analytic functions of $h$.
 Using lemma \ref {lem:ana} we conclude that
all the coefficients of $ J_{q-1} f_{ (0,0)} $ are analytic in $ h $.

If follows that $ \tau_m (h) $ is analytic in $ h $, leading immediately to
the next corollary.
 \qed

\begin {corollary}
On each $ I_j^q $ and for $1 \le m \le \left [ \frac{q}{2} \right ] -1 $,
the set $ \{ h \, / \, \tau_m (h) = 0 \} $ is either the entire $ I_j^q$
or a discrete set.
 \end {corollary}

In order to prove the existence of islands we use the natural recurrence on the order of the ressonaces.

We begin by analysing the zeros of $\tau_1$ on the non resonat intervals $I_j^4$. 
If $\tau_1 = 0$ only on a discrete subset of each $I_j^4$, $Pan(n_0,a_0,h)$
has elliptic island except for a discrete set of values of $h$ (which can
be smaller than the union of the discrete subsets of zeros of $\tau_1$ and the
values of resonance up to order 4, since on the discrete subsets a non
resonnat value of higher order may have a non zero Birkhoff coefficient).

If $\tau_1$ is identically zero on one of those intervals, we proceed to the next step, 
aplying the same analysis to the zeros of $\tau_2$.

We continue the recurrence and it will end up in a finite number of steps if
for some order of resonance the last Birkhoff coefficient does not vanish
identically on a whole non resonant interval. Otherwise, all the Birkhoff
coefficients will vanish on at least one open interval, bounded by resonant
values of $h$.
In this last case, since $\frac{\partial \Delta_n} {\partial h} > 0$, the rotation
number $\rho (h)$ of $Pan (n_0,a_0,h)$ is not constant and there exists $h_0$ such
that $\rho (h_0)$ is diophantine. As $f$ is analytic, it is conjugate to a
rotation  and there will be invariant curves \cite{kn:her}.

We conclude that:
\begin {theorem}
Given $ n $ and $ a > \alpha_n^1 $, there are at least countably many values of
$ h $ in $ \cup I_j^4 $ such that $ Pan ( n , a , h ) $ has an elliptic
island.   
\end {theorem}

\Remark As Moeckel proved in \cite{kn:moe}, in a generic one-parameter family
of area preserving maps with elliptic fixed points, the first Birkhoff
coefficient $\tau_1$ varies from $-\infty$ to $+\infty$ as the rotation number varies 
from 0 to 1/3 or from 2/3 to 1. So, we do not expect $\tau_1(h)$ to be always
different from zero, neither to vanish identically. Furthermore, generically,
at a zero of $\tau_1$, a higher Birkhoff coefficient will not vanish.
So, although we can handle the case $ \tau_s (h_0) = 0 $, $ \forall s $, we do
not expect it to happen in our case.

 \section {Bounds for the existence of islands}
 \subsection {The case $ a< \sqrt 2 $ }

As shown in proposition \ref {prop:region} and in \cite {kn:cana}, if $ (a, h) $
is in the region of ellipticity of $ Pan (n, a, h) $ then 
 $ a \geq \alpha_n^1 = \sqrt { \frac {2n + 2} {n + 2}} $. 
So, given $ a \in (1, \sqrt 2) $, there is only a finite number of periods
 $ 4 + 2n $ such that $ Pan (n, a, h) $ can be elliptic.
 More precisely, $ n \leq \frac {2 (a^2-1)} {2-a^2} $. 
 
Let $H(a)$ be the maximum of $h^1_n(a)$ 
for those periods. As proved in \cite {kn:cana}, $H(a)$ is a lower bound for chaos.
By theorem~1, it is also an upper bound for the existence of elliptic islands for the
Pantographic family.

 \subsection {The case $ a> \sqrt 2 $ }

On the other hand, if $ a> \sqrt 2 $, $ Pan (n, a, h) $ can be elliptic for any
period $n$. Moreover, for each $ n $ and $q$,
$ \cup I_j^q \subset (h_n^0 (a),h_n^1 (a)) $, with 
 $ h_n^0 (a) = n a \sqrt {a^2-2} $. We also have that 
 $ (h_n^0 (a), h_n^1 (a)) \cap (h_n^0 (a), h_{n + 1}^0 (a)) $  
is a non empty open interval. So we can find $ h $,
 $ n a \sqrt {a^2-2} < h < (n + 1)a \sqrt {a^2-2} $ such that
 $ Pan (n, a, h) $ has an elliptic island. This proves the following

\begin {theorem}
Given $ a> \sqrt 2 $ there is no upper bound on $ h $ for the existence of
elliptic islands on the elliptical stadium billiard.
 \end {theorem}

However, as can be seen in figure \ref {fig:gaps}, for values of $ a $ further away
from $\sqrt{2}$, the strips of ellipticity are disjoint. In these gaps all
pantographic orbits are hyperbolic, having, thus, no islands. 

\bfig[h]
\putfigurewithtex gaps.ps  gaps.txt 8 16
 \caption {Gaps between the strips of ellipticity}
 \label {fig:gaps}
 \efig

Nevertheless, in our simulations   other islands appear, obviously
corresponding to different periodic orbits. In figure \ref{fig:malucas} we
exemplify this fact for $a=2$ and $h=2$, a value located in the gap between
the strips of ellipticity for $n=0$ and $n=1$ (see figure~\ref {fig:gaps}). We
show three non-pantographic orbits and their islands and the whole phase
space, where we can see many other islands surrounded by what seems to be a
chaotic sea. As far as our results indicate and our simulations show,
this should be the typical picture for the phase space when $ a> \sqrt 2 $.

\bfig
\vskip 7.0 truecm 
\includegraphics{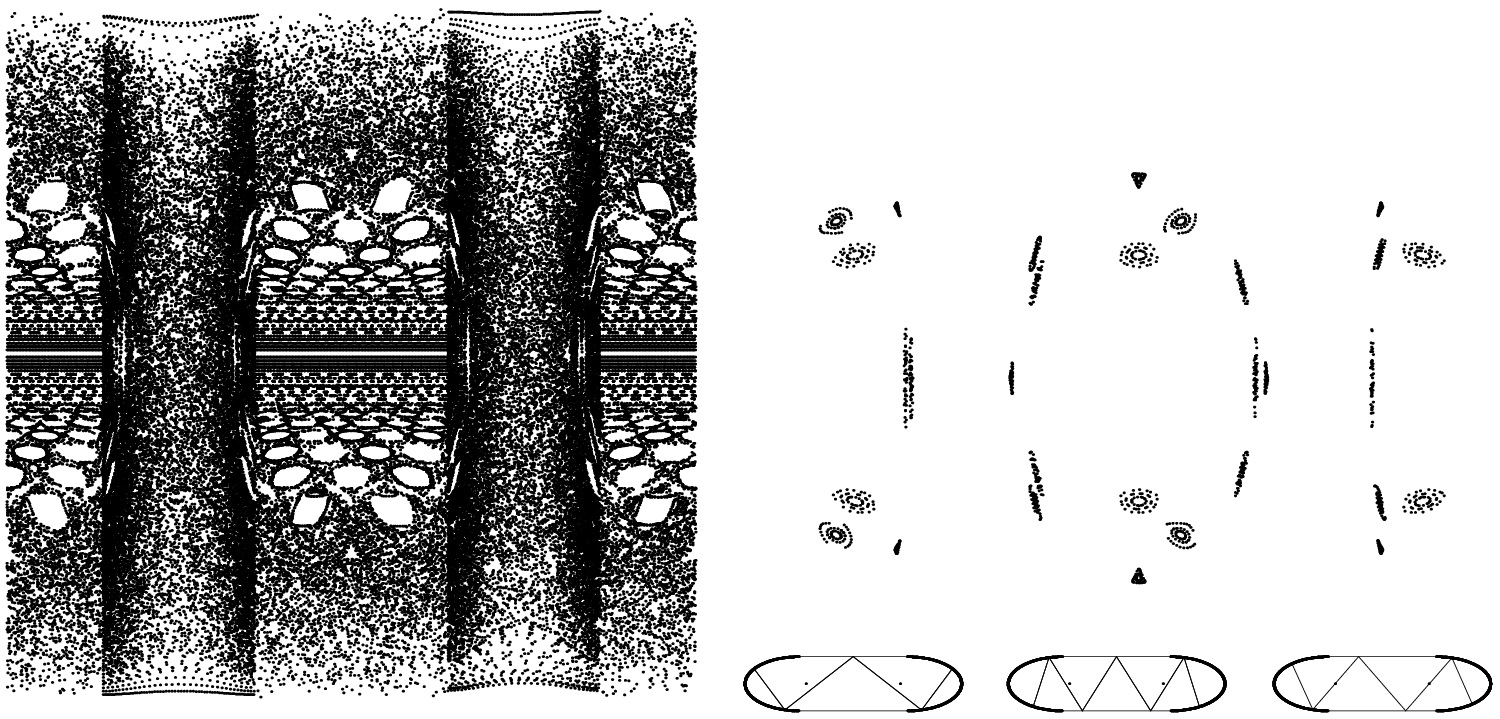}
 \caption {Phase space for $(a,h)=(2,2)$}
 \label {fig:malucas}
 \efig

\section {Appendix}

\subsection{Some properties of $t(n,a,h)$}

We called $ U_n $ the open region in the parameter plane where
 $ Pan (n, a, h) $ exists:
$$ U_0 = U_1 = \{ (a, h)/ \, a>1, h>0 \}$$ 
$$ U_n = \{ (a, h)/ \,1<a \leq 2,h>0 \} \cup \{ (a, h)/ \, a>2 ,
h> (n-1) \sqrt {a (a-2)} \},$$ 
for $ n \geq 2 $.

The following lemma gives some useful information about  
 $ t (n, a, h) $, the solution of (\ref {eq:pan}) or (\ref {eq:y}).

\begin {lemma}
Given $ n \geq 0 $, let $ (a, h) \in U_n $ and $ t (n, a, h) $ be the unique
solution of (\ref {eq:y}). Then 
\begin {itemize}
\item $ \forall n $, when $ h \to + \infty $, $ t (n, a, h) \to \frac {1} {a} $.
\item For $ n = 0 $ 
\begin {itemize}
\item if $ 1<a \leq \sqrt 2 $, when $ h \to 0^ + $, $ t (0, a, h) \to + \infty $ 
\item if $ a> \sqrt 2 $, when 
$ h \to 0^ + $, $ t (0, a, h) \to \frac {1} { \sqrt {a^2-2}} $ 
 \end {itemize} 
\item For $ n = 1 $ 
\begin {itemize}
\item
if $ 1<a \leq 2 $, when $ h \to 0^ + $, $ t (1, a, h) \to + \infty $ 
\item
if $ a>2 $, when $ h \to 0^ + $, 
$ t (1, a, h) \to \frac {1} { \sqrt {a (a-2)}} $ 
\end {itemize}
\item For $ n \geq 2 $ 
\begin {itemize}
\item if $ 1<a \leq 2 $, when $ h \to 0^ + $, $ t (n, a, h) \to + \infty $ 
\item if $ a>2 $, when $ h \to { (n-1) \sqrt {a (a-2)}}^ + $, $ t (n, a, h) \to
  \frac {1} { \sqrt {a (a-2)}} $ 
 \end {itemize}
 \end {itemize}
 \end {lemma}

\Proof
Since $ \frac { \partial t} { \partial h}<0 $ and $ \frac {1} {a} < t $,
$ {\displaystyle \lim_{h \to + \infty } t } $ exists.
From equation (\ref {eq:pan}),
\be 
a t - \frac {1} {at} = 
2 \frac {n \sqrt {1 + t^2} + t} {h \sqrt {1 + t^2} + a} \to 0
\en 
as $ h \to + \infty $ and so $ t \to \frac {1} {a} $.

To study the limit as $ h \to 0^ + $, let us take $ x = \frac {1} {t} $.
Equation 
 (\ref {eq:pan}) becomes
$$
x (2n \sqrt {1 + x^2} + x^2- (a^2-2)) = 
(a^2 - x^2) \sqrt {1 + x^2} \, \frac {h} {a}.$$ 
Since $ 0<x< \frac {1} {a} $, $ x (2 n \sqrt {1 + x^2} + x^2- (a^2-2))>0 $.

Let $ \overline x = {\displaystyle \lim_{h \to 0^ + } } x $.
We have that $ 0 \leq \overline x \leq a $ and 
$ 
\overline x (2n \sqrt {1 + { \overline x}^2} + { \overline x}^2- (a^2-2)) = 0 $.

For $ n = 0 $, if $ 1<a \leq \sqrt 2 $, $ \overline x = 0 $ is the unique
solution of this equation and
$ {\displaystyle \lim_{h \to 0^ + } }t (0, a, h) = + \infty $.
If $ a> \sqrt 2 $ we have a new
solution $ \overline x = \sqrt {a^2-2} $. But for $ 0<x< \sqrt {a^2-2} $,
 $ x (x^2- (a^2-2))<0 $.
 So $  {\displaystyle \lim_{h \to 0^+ } } x = \sqrt {a^2-2} $ and 
 $ {\displaystyle \lim_{h \to 0^+ } } t (0, a, h) = 
 \frac {1} { \sqrt {a^2-2}} $.

For $ n = 1 $ we have $ \overline x(2 \sqrt {1 + { \overline x}^2} + { \overline
x}^2- (a^2-2)) = 0 $. If $ a^2-2 \leq 2 $, i.e. $ a \leq 2 $,
the unique solution is 
 $ \overline x = 0 $ and
 $ {\displaystyle\lim_{h \to 0^+ }} t (1, a, h) = + \infty $.
 If $ a>2 $ the second
solution is $ \overline x = \sqrt {a (a-2)} $.
As above, if $ 0 < x < \sqrt {a (a-2)} $,
$ x(2 \sqrt {1 + x^2} + x^2- (a^2-2))<0 $ and so
${\displaystyle \lim_{h \to 0^+ }} t (1,a, h) = 
 \frac {1} { \sqrt {a (a-2)}} $.

For $ n \geq 2 $, we remark first that if $ k > l $ then
$ t (k, a, h) > t (l, a, h) $. So, for
 $ 1<a \leq 2$ $ {\displaystyle \lim_{h \to 0^+ }} t (n, a, h) = + \infty $.
 For $ a>2 $, the limit as 
 $ h \to (n-1) \sqrt {a (a-2)} $ is the unique solution of 
 equation~(\ref {eq:pan})
for $ h = (n-1) \sqrt {a (a-2)} $ which is 
$ t = \frac {1} { \sqrt {a (a-2)}} $. \qed

{\footnotesize Remark: 
When $ h \to 0^ + $, the elliptical stadium becomes an ellipse.
$ t \to + \infty $ 
means that the pantographic orbit goes to the elliptic periodic orbit which
corresponds to the minor axis of the ellipse.\\
Let us call { \it pantographic-like} orbits in the elliptical billiard the
periodic trajectories that have vertical segments both at left and right
extremes. As can be seen in \cite {kn:cmp}, the 4-periodic pantographic-like
orbit exists if $ a> \sqrt 2 $ and the 6-periodic if $ a>2 $, they are parabolic
and their position is given, respectively, by $ t = \frac {1} { \sqrt {a^2-2}} $
and $ t = \frac {1} { \sqrt {a (a-2)}} $.
They are the calculated limits of the 4 and
6-periodic pantographic orbits of the elliptical stadium.}

 \subsection {Proof of the lemma 2}

For a fixed $ n $, let $ (a, h) \in U_n $ and $ \lambda_n (a, h) $ be the
solution of~(\ref {eq:pan}), $ \beta $ be the angle, with the normal, of the
outgoing trajectory at $ P = (a \cos \lambda_n + h, \sin \lambda_n) $ and $ s $
the corresponding arc length. 

 Let 
$$ 
 \Delta_n (a, h) = \left ({l_1 \, K \over \cos \beta} -1 \right ) \, 
 \left ({l_2 \, K \over \cos \beta} -1 \right )
$$ 
 where 
 $ l_1 = 2 \sin \lambda_n $, 
 $ l_2 = 2 \sqrt { (h + a \cos \lambda_n)^2 + (n + \sin \lambda_n)^2} $ and 
 $ K = a/ (a^2 \sin^2 \lambda_n + \cos^2 \lambda_n)^{3/2} $. 

Let $ \delta_1 (a, h) = \left ({l_1 \, K \over \cos \beta} -1 \right ) $ and
 $ \delta_2 (a, h) = \left ({l_2 \, K \over \cos \beta} -1 \right ) $.

\begin {lemma} For every $ n \geq 0 $ the function $ \delta_1 (a, h) $ has the
following properties:
\begin {enumerate}
\item $ \delta_1 (a, h)>0 $ for $ 1<a< \sqrt 2 $ and $ h>0 $.
\item
$ {\displaystyle \lim_{h \to 0^ + } } \delta_1 (a, h) = \frac {2} {a^2}-1 $ for
$ 1<a< \sqrt 2 $.
\item $ \delta_1 (a, n a \sqrt {a^2-2}) = 0 $ for $ a \geq \sqrt 2 $.
\item $ \frac { \partial \delta_1} { \partial h}>0 $ for $ (a, h) \in U_n $.
 \end {enumerate}
 \label {lem:d1}
 \end {lemma}
 \Proof
We have 
 $ \delta_1 (a, h) = {l_1 \, K \over \cos \beta}-1 =
 2 {1 + t^2 \over 1 + a^2 t^2}-1 $ 
and properties 1 and 2 follow immediately.

If $ a \geq \sqrt 2 $, $ \delta_1 = 0 $ implies $ (a^2-2) t^2-1 = 0 $ and
 $ t = \frac {1} { \sqrt {a^2-2}} $. From equation (\ref {eq:pan})
 $ h = n a \sqrt {a^2-2} $ and property 3 follows.

Since 
$ \frac { \partial \delta_1} { \partial h} = 
\frac { \partial \delta_1} { \partial t}
 \frac { \partial t} { \partial h} $, as 
 $ \frac { \partial \delta_1} { \partial t} = 
 \frac {4t (1-a^2)} { (1 + a^2 t^2)^2}<0 $, for $ a>1 $, and
 $ \frac { \partial t} { \partial h}<0 $, 
 $ \frac { \partial \delta_1} { \partial h}>0 $. \qed

\begin {lemma} For every $ n \geq 0 $ the function $ \delta_2 (a, h) $ has the
following properties:
\begin {enumerate}
\item $ \delta_2 (a, h)> \delta_1 (a, h)>0 $ for $ (a, h) \in U_n $.
\item
$ {\displaystyle \lim_{h \to 0^ + }}
 \delta_2 (a, h) = 2 \frac {n + 1} {a^2}-1 $ for $ 1<a< \sqrt 2 $.
\item
${\displaystyle \lim_{h \to + \infty}} \delta_2 (a, h) = + \infty $ for $ 1<a $.
\item $ \frac { \partial \delta_2} { \partial h}>0 $ for $ (a, h) \in U_n $.
 \end {enumerate}
 \label {lem:d2}
 \end {lemma} 
 \Proof
Since $ l_2 > l_1 $, $ \delta_2> \delta_1 $.

For $ 1<a< \sqrt 2 $, when $ h \to 0 $, $ t (n, a, h) \to + \infty $.
So $ l_2 \to 2 (n + 1),
K \to 1/a^2 $ and $ \cos \beta \to 1 $ and property 2 follows.

Property~3 is obvious since $ l_2 \to \infty $ as $ h \to \infty $ and all the 
other quantities are bounded. 

By definition 
$ \delta_2 (a, h) = 
l_2 (a, h, t (a, h)) \frac {K (t (a, h))} { \cos \beta (t (a, h))}-1 $
and $ \frac { \partial \delta_2} { \partial h} = 
 \frac {K} { \cos \beta} \frac { \partial l_2} { \partial h} +
 \frac { \partial} { \partial t}
 \left ({l_2 \, K \over \cos \beta} \right) \frac { \partial t} { \partial h}. $ 

We have that $ \frac { \partial t} { \partial h}<0 $ and 
 $ \frac { \partial t} { \partial \lambda} > 0 $. The curvature $ K>0 $ and for 
 $ 0 < \lambda < \pi / 2 $, $ \frac { \partial K} { \partial \lambda} < 0 $. So
 $ \frac { \partial K} { \partial t}<0 $. As $ 0< \beta < \pi /4 $, 
 $ \cos \beta > 0 $.
As $ \tan \beta = 1/at $, $ \frac { \partial} { \partial t} \cos \beta > 0 $.
This implies that 
 $ \frac { \partial} { \partial t} \left (\frac {K} { \cos \beta} \right) < 0 $.

We have that 
 $ \frac {1} {4}l_2^2 = (h + a \cos \lambda)^2 + (n + \sin \lambda)^2 $. So
 $ \frac {1} {2}l_2 \frac { \partial l_2} { \partial h} = 
 2 (h + a \cos \lambda)>0 $ in
 $ U_n $, implying that $ \frac { \partial l_2} { \partial h}>0 $.
We also have that $ \frac {1} {4} \frac { \partial l_2^2} { \partial \lambda} = 
-2a \sin \lambda (n + \sin \lambda) \, (\tan 2 \beta - \tan \beta)<0 $ 
 as $ 0< \beta < \pi /4 $. So $ \frac { \partial l_2} { \partial \lambda} = 
 \frac { \partial l_2} { \partial t} 
 \frac { \partial t} { \partial \lambda} < 0 $, and 
 $ \frac { \partial l_2} { \partial t}<0 $.

This shows that $ {l_2 \, K \over \cos \beta} $ is the product of two positive
decreasing functions of $ t $ and so $ \frac { \partial} { \partial t}
 \left ({l_2 \, K \over \cos \beta} \right) < 0 $.

We conclude that $ \frac { \partial \delta_2} { \partial h}>0 $. \qed

We have defined $ \tilde U_n = \{ (a, h)/ \,1<a< \sqrt 2,h>0 \} \cup
 \{ (a, h)/ \, a \geq \sqrt 2, h > n a \sqrt {a^2-2} \} \subset U_n $.
 \\
 \\
 { \bf Lemma 2}
 { \em For every $ n \geq 0 $ the function $ \Delta_n (a, h) $ has the following
properties:\begin {enumerate} 
\item $ \Delta_n (a, h) | _{ \tilde U_n } > 0 $ 
\item $ \frac { \partial \Delta_n} { \partial h} | _{ \tilde U_n } > 0 $ 
\item $ {\displaystyle \lim_{h \rightarrow + \infty}}
 \Delta_n (a, h) = + \infty $ 
\item for $ 1<a< \sqrt 2 $,
$ {\displaystyle \lim_{h \rightarrow 0} } \Delta_n (a, h) = 
  L_n (a) = (\frac {2} {a^2}-1) (\frac {2 (n + 1)} {a^2}-1)> 0 $ \\
for $ a \geq \sqrt 2 $, 
$ {\displaystyle \lim_{h \rightarrow n a \sqrt {a^2-2}} }
\Delta_n (a, h) = 0 $.
 \end {enumerate} 
}
 \Proof
In $ \tilde U_n $, $ \Delta_n $ is the product of two positive increasing
functions of $ h $ and properties 1 and 2 follows.
Properties 3 and 4 follow immediately
from lemmas \ref {lem:d1} and \ref {lem:d2}. \qed

\ack
We want to thank A.Chenciner, M.Herman, C.Grotta Ragazzo and R.Roussarie. We
also thank the support of CNPq during this work.

email: syok@mat.ufmg.br and sonia@mat.ufmg.br 


\begin {thebibliography} {99} 
 \bibitem {kn:bun} 
 L.~A.~Bunimovich: Conditions on stochasticity of two-dimensional billiards. 
 Chaos { \bf 1}, 187-193, (1991) 
 \bibitem {kn:cana}
E.~Canale, R.~Markarian, S.~Oliffson Kamphorst, S.~Pinto de Carvalho: A lower
bound for chaos on the elliptical stadium. Physica D { \bf 115}, 189-202 (1998)
 \bibitem {kn:don} 
V.~J.~Donnay: 
Using integrability to produce chaos: billiards with positive 
entropy. Comm. Math. Phys. { \bf 141}, 225-257 (1991)
\bibitem {kn:her}
M. Herman: private communication
\bibitem{kn:cmp} 
R.~Markarian, S.~Oliffson Kamphorst, S.~Pinto de Carvalho: 
Chaotic Properties 
of the Elliptical Stadium. Comm. Math. Phys. {\bf 174}, 661-679 (1996) 
\bibitem {kn:mar} 
R.~Markarian: Introduction to the ergodic theory of plane billiards. 
In: { \em Dynamical Systems}, Santiago de Chile, 1990.  
 \bibitem {kn:moe}
R.~Moeckel: Generic bifurcations of the twist coefficient, Erg.Th.Dyn.Syst. 
 { \bf 10}, 185-195 (1990)
 \bibitem {kn:mos}
J.~Moser: { \em Stable and random motions in dynamical systems}, PUP, Princeton,
1973.
  \end {thebibliography}
 \end {document}